\documentclass[a4]{article}

\usepackage{amsfonts}
\usepackage{bm}
\usepackage{bbm}
\usepackage{graphicx}

\newcommand{\btheta}{\boldsymbol{\theta}}
\newcommand{\blambda}{\boldsymbol{\lambda}}
\newcommand{\bmu}{\boldsymbol{\mu}}

\newcommand{\bl}{\boldsymbol{l}}

\newcommand{\bp}{\boldsymbol{p}}
\newcommand{\bq}{\boldsymbol{q}}

\newcommand{\bI}{\boldsymbol{I}}
\newcommand{\bx}{\boldsymbol{x}}
\newcommand{\bX}{\boldsymbol{X}}
\newcommand{\by}{\boldsymbol{y}}

\newcommand{\R}{{\mathbb R}}
\newcommand{\Z}{{\mathbb Z}}
\newcommand{\N}{{\mathbb N}}

\begin{document}

\title{Persistence time of SIS infections in heterogeneous populations and networks}
\author{Damian Clancy\\
Department of Actuarial Mathematics and Statistics\\ Maxwell Institute for Mathematical Sciences\\ Heriot-Watt University\\ Edinburgh\\ EH14 4AS\\ UK\\
d.clancy@hw.ac.uk}
\maketitle

\begin{abstract}
For a susceptible-infectious-susceptible (SIS) infection model in a heterogeneous population, we present simple formulae giving the leading-order asymptotic (large population) behaviour of the mean persistence time, from an endemic state to extinction of infection.
Our model may be interpreted as describing an infection spreading through either (i)~a population with heterogeneity in individuals' susceptibility and/or infectiousness; or (ii)~a heterogeneous directed network.
Using our asymptotic formulae, we show that such heterogeneity can only reduce (to leading order) the mean persistence time compared to a corresponding homogeneous population, and that the greater the degree of heterogeneity, the more quickly infection will die out.
\vspace*{4mm}

{\bf Keywords: Stochastic epidemic models; large deviations; endemic fade-out; directed configuration model; superspreaders}
\end{abstract}

\section{Introduction}
In modelling endemic infections, a quantity of particular interest is the persistence time until infection dies out from the population.
For discrete state-space Markov chain models, the expected persistence time for an infection that has become endemic in the population (i.e. starting from quasi-stationarity) may be found as an eigenvalue of the transition rate matrix.
However, for large populations and for more complicated models, numerical computation of  this exact solution can be very time-consuming, and may also suffer from numerical instability.
Moreover, it is not straightforward to use this eigenvalue characterization to investigate, for instance, the effect of population heterogeneities upon the expected persistence time.
Approximation methods are therefore required.
For a number of infection models, it has been shown that, denoting by $N$ the typical size of the population, the expected time from endemicity to extinction, $\tau$, is asymptotically given by an expression of the form
\begin{eqnarray}
\tau &\sim& {C \over \sqrt{N}} \exp ( AN )
\label{general_asymptotic_formula}
\end{eqnarray}
where the values of $A,C$ depend upon parameters of the model, but not upon~$N$.
We assume here that the process is super-critical, so that long-term endemicity is possible.

For the classic susceptible-infectious-susceptible (SIS) model of Weiss and Dishon~\cite{WD71}, Andersson and Djehiche~\cite{AD98} found simple explicit expressions for both $A$ and $C$ in terms of the basic reproduction number $R_0$ (the expected number of secondary cases caused by a typical primary case in an otherwise susceptible
population), under the assumption of super-criticality (that is, $R_0 > 1$); specifically, $A = (1/R_0) - 1 + \ln R_0$ and $C = R_0 \sqrt{2\pi} / ( R_0 - 1 )^2$, assuming time is scaled such that individual infectious periods are of mean~1.
This was extended by Ball,~Britton and~Neal~\cite{BBN16} to allow for a general infectious period distribution in place of the exponential distribution assumed by~\cite{AD98}; they showed that leading-order behaviour is unchanged, so that $A = (1/R_0) - 1 + \ln R_0$ as before, while the value of $C$ depends upon the infectious period distribution, and may be straightforwardly evaluated provided this distribution is known.
Pre-dating the above work, van Herwaarden and Grasman~\cite{VHG95} showed that relationship~(\ref{general_asymptotic_formula}) holds true for a particular suceptible-infectious-removed (SIR) infection model.
In this case, however, evaluation of the constant~$A$ requires numerical solution of a system of ordinary differential equations, while no method for evaluating~$C$ is given.

The system of ordinary differential equations used in~\cite{VHG95} to evaluate~$A$ may be regarded as the equations of motion corresponding to a particular Hamiltonian system.
More recently, a number of authors~\cite{AM10,AM16,DMRH94,EK04,KM08,LSS14} have applied this Hamiltonian approach to a range of infection models to derive results of the form
\begin{eqnarray}
\lim_{N \to \infty} {\ln \tau \over N} &=& A \label{general_rough_result} .
\end{eqnarray}
Equation~(\ref{general_rough_result}) is not as precise as  relationship~(\ref{general_asymptotic_formula}),
but does at least give the leading-order behaviour of $\tau$ in the large population ($N \to \infty$) limit. 
Evaluation of~$A$ generally requires numerical solution of the equations of motion, and consequently much of the research effort has focused upon developing efficient numerical procedures.

We shall apply the Hamiltonian approach to approximate the expected persistence time from endemicity, $\tau$, for an SIS model incorporating heterogeneity in individual infectivities and susceptibilities. 
Such heterogeneity is a common feature of real-world infections.
For instance, for a number of infections (eg SARS) it has been hypothesised that  there exists a subgroup of `super-spreaders' within the population, being individuals of higher infectivity than the rest.
Heterogeneous susceptibilities may arise, for instance, through individuals having differing histories of prior exposure to infection or vaccination.
Alternatively, our model may  be interpreted as a model for infection spreading on an uncorrelated (that is, with no correlations between degrees of neighbouring individuals) directed network~\cite{DGM08}.

In contrast to almost all previous work, we are able to find an explicit formula for the constant~$A$ in equation~(\ref{general_rough_result}), at least provided the heterogeneity is in either infectivity, or susceptibility, but not both.
As well as being much quicker and easier to evaluate than the solution to a (typically high-dimensional) system of ordinary differential equations, a further advantage of such an explicit formula is that it may be used to establish qualitative results about the effects of model assumptions.
Specifically, we investigate the effect of increasing heterogeneity upon the persistence time of infection in the population.

The remainder of the paper is structured as follows.
In section~2, we define precisely our heterogeneous population SIS model, and describe how it may be interpreted as approximating a directed network model.
In section~3 we recall some general theory that will be required in the sequel.
Our main result, theorem~1, is derived in section~4, establishing explicit asymptotic formulae for $\ln \tau$ in the large-population limit, provided that heterogeneity is in either infectivity or susceptibility, but not both.
Using these explicit formulae, we go on in section~5 to demonstrate that the greater the level of heterogeneity in either infectivity or susceptibility, the more rapidly extinction of infection will occur (on average, to leading order in a large population).
In the case that both heterogeneities are present simultaneously, numerical work (figure~\ref{contour}) suggests that mean persistence time is again maximised in the homogeneous-population case.
In section~6 we demonstrate that if heterogeneity is in susceptibilities, our asymptotic formula for $\ln \tau$, and hence also our conclusion that greater heterogeneity reduces  mean persistence time, remain valid for more general infectious period distributions than the classic exponential distribution.
We present numerical evidence (figure~\ref{const_inf_periods}) suggesting that this is also true when instead heterogeneity is in infectivities.
Finally, in section~\ref{discussion}, we discuss the directed network interpretation of our results (theorem~4) and suggest some directions for further work.

\section{The SIS infection model in a heterogeneous population or directed network}
We first formulate our model in terms of a population divided into a fixed number of groups, and then describe how the same model may be interpreted as modelling an infection spreading on a directed network.

Consider a closed population of $N$ individuals divided into $k$ groups, with group~$i$ ($i=1,2,\ldots,k$) consisting of $N_i$ individuals.
Denote by $f_i = N_i / N$ the proportion of the population belonging to group $i$, so that $\sum_i f_i = 1$.
When a group~$i$ individual becomes infected, it remains so for a time distributed as an exponential random variable with mean $1/\gamma$  (assumed for simplicity to be the same for each group).
During this infectious period, the group~$i$ infective makes contacts with each individual in each group~$j = 1 , 2 , \ldots , k$ at the points of a Poisson process of rate $\beta \lambda_i \mu_j / N$, where $\beta$ is some overall measure of infectiousness, $\lambda_i$ represents the infectivity of group~$i$ individuals, and $\mu_j$ represents the susceptibility of group~$j$ individuals. 
(The assumption that the group~$i$ to group~$j$ infection rate factorises in this way is sometimes referred to as `separable mixing'.)
Without loss of generality, we scale the $\lambda_i$, $\mu_j$ values so that $\sum_i \lambda_i f_i = \sum_j \mu_j f_j = 1$.
These Poisson processes and infectious periods are all mutually independent.
If a contacted individual is susceptible, then it becomes infected (and infectious); if the contacted individual is already infected then the contact has no effect.
Thus the process $\left\{ \bI (t) = \left( I_1 (t) , I_2(t) , \ldots , I_k (t) \right) : t \ge 0 \right\}$ is a continuous-time Markov chain with transition rates given in table~\ref{k_group_SIS}, and the number of susceptible individuals in group~$j$ at time~$t \ge 0$ is $S_j (t) = N_j - I_j (t)$.
We will assume throughout that $\beta , \gamma > 0$, and that $f_i , \lambda_i , \mu_i > 0$ for all $i$.
Note that our model is a special case of the model of~\cite{CP13}, although we use slightly different notation here.
The basic reproduction number $R_0$ is given by the  dominant eigenvalue of the matrix $M$ with entries $m_{ij} = \beta  \lambda_i \mu_j f_j$, so that
\begin{eqnarray*}
R_0 &=& {\beta \over \gamma} \sum_{i=1}^k \lambda_i \mu_i f_i .
\end{eqnarray*}

\renewcommand{\arraystretch}{2}
\begin{table}
\begin{tabular}{lll} \hline
Event & State transition & Transition rate \\ \hline
Infection in group $j$ & $I_j \to I_j + 1$ & ${\beta \over N} \left( \sum_{m=1}^k \lambda_m I_m \right) \mu_j (N_j - I_j)$ \\
Recovery in group $j$ & $I_j \to I_j - 1$ & $\gamma I_j$ \\
\hline
\end{tabular}
\caption{Transition rates for the $k$-group SIS model.}\label{k_group_SIS}
\end{table}
\renewcommand{\arraystretch}{1}

We now describe how the above model may be interpreted as describing infection spreading through a network.
Each of the $N$ individuals in the population is assigned an in-degree~$d_{\mbox{in}}$ and out-degree~$d_{\mbox{out}}$ according to some joint probability mass function $p \left(d_{\mbox{in}},d_{\mbox{out}} \right)$ on $\Z_+^2$.
These degrees are assigned independently to distinct individuals, but the in and out degree need not be independent for a single individual.
To each individual we attach `stubs', or half-edges, with $d_{\mbox{in}}$ stubs pointing inwards and $d_{\mbox{out}}$ stubs pointing outwards.
Inward-pointing stubs are then paired with outward-pointing stubs throughout the population, to create links between individuals.
In order that this process can produce a valid network, with no left over half-edges, we clearly require that $E \left[ d_{\mbox{in}} \right] = E \left[ d_{\mbox{out}} \right]$.
We do not concern ourselves with the precise mechanism by which stubs are paired off (see~\cite{BJM07,CO13} for relevant discussion); rather, we shall simply assume that the resulting network is uncorrelated, so that the so-called `annealed' network approximation is valid for an ensemble of such networks.
This is a mean-field approximation for heterogeneous networks, and may be described as follows
(see~\cite{DGM08} for more comprehensive discussion).

Denote by $\kappa$ the rate at which infection transmits along each link from an infectious individual to a susceptible individual.
Suppose for simplicity that there are a finite number~$k$ of $\left( d_{\mbox{in}} , d_{\mbox{out}} \right)$ pairs having non-zero probability, and define a bijective function $c \left( d_{\mbox{in}} , d_{\mbox{out}} \right) : \Z_+^2 \to \left\{ 1,2,\ldots,k \right\}$ that assigns a unique number to each of the possible $\left( d_{\mbox{in}} , d_{\mbox{out}} \right)$ pairs.
We say an individual is of `group~$j$' if they have degrees $\left( d_{\mbox{in}} , d_{\mbox{out}} \right) = c^{-1} (j)$, and define $d_{\mbox{in}} (j) , d_{\mbox{out}} (j)$ to be the in and out degrees, respectively, of a group~$j$ individual.
For $j = 1,2,\ldots,k$, denote by $N_j$ the total number of group~$j$ individuals in the population, and by $I_j (t)$ the number of group~$j$ individuals who are infectious at time~$t$.
Then under the annealed network approximation, the total rate at which group~$j$ individuals become infected is given by 
\begin{eqnarray*}
{\kappa \over N E \left[ d_{\mbox{in}} \right]} \left( \sum_{m=1}^k d_{\mbox{out}} (m) I_m  \right)  d_{\mbox{in}} (j) \left( N_j - I_j \right) .
\end{eqnarray*}
When individuals become infected, they remain so for an exponentially distributed time of mean~$1/\gamma$ before returning to the susceptible state.

This network model may be approximated by the $k$-group model with transition rates given in table~\ref{k_group_SIS} by taking
\begin{eqnarray*}
\beta &=& \kappa E \left[ d_{\mbox{out}} \right] , \\
f_j &=& p \left( c^{-1} (j) \right) , \\
\mu_j &=& d_{\mbox{in}} (j) / E \left[ d_{\mbox{in}} \right] , \\
\lambda_j &=& d_{\mbox{out}} (j) / E \left[ d_{\mbox{out}} \right] ,
\end{eqnarray*}
for $j=1,2,\ldots,k$.

The undirected version of the above annealed network approximation (with $\blambda = \bmu$) has been studied by~\cite{HS16}, via numerical solution of Hamilton's equations of motion, equations~(\ref{eqn_motion_y},\ref{eqn_motion_theta}) below.

In the next section we present some relevant general theory, before going on in section~4 to apply these general methods to the model described above.

\section{General theory regarding persistence time from endemicity}
\label{theory}
Consider an infection modelled by a continuous-time Markov process $\left\{ \bX (t) : t \ge 0 \right\}$ on finite state-space~$S \subset \Z^k$ with transition rate matrix~$Q$.
Suppose that~$S$ is made up of an absorbing set of states~$A$ (corresponding to absence of disease) and a single transient communicating class~$C$.
We denote by $Q_C$ the transition rate matrix restricted to~$C$.
Then the infection will almost surely die out (i.e. the process will leave~$C$) within finite time, and~\cite{DS67}
there exists a unique quasi-stationary distribution $\bq = \left\{ q_{\bx} : \bx \in C \right\}$ such that, for any initial state within~$C$,
\begin{eqnarray*}
q_{\bx} &=& \lim_{t \to \infty} \Pr \left( \bX(t) = \bx \left| \bX (t) \in C \right. \right)
\mbox{ for } \bx \in C .
\end{eqnarray*}
That is, provided the infection does not die out, it will settle to the endemic  distribution~$\bq$.
The distribution $\bq$ may be found as the unique solution of
\begin{eqnarray}
\bq Q_C &=& - (1/\tau) \bq \mbox{ with } \sum_{\bx \in C} q_{\bx} = 1 , \label{QSD}
\end{eqnarray}
where $-(1/\tau)$ is the eigenvalue of $Q_C$ with largest real part.
The time to extinction from quasi-stationarity is exponentially distributed with mean~$\tau$.

Although $\tau$ may be computed exactly from equation~(\ref{QSD}), this can become impractical when the state-space is large, and it is not straightforward from~(\ref{QSD}) to establish qualitative results.
Approximation methods are therefore valuable, and in particular, methods from Hamiltonian statistical mechanics may be used to study the leading order asymptotic (large population) behaviour of $\tau$, as follows.

Suppose that $\bX(t)$ is a density-dependent process in the sense of chapter~11 of~\cite{EK05}; that is, the transition rates are of the form
\begin{eqnarray}
P \left( \bX ( t + \delta t ) = \bx + \bl \mid \bX (t) = \bx \right) &=& N W_{\bl} \left( {\bx \over N} \right) + o ( \delta t ) \mbox{ for } \bx \in S,\ \bl \in L , \label{rates}
\end{eqnarray} 
for some functions $W_{\bl} : \R^k \to \R^+$, where $L$ is the set of possible jumps from each state $\bx \in S$ and $N$ is some parameter indicating overall size of the system (in our applications, $N$ will be the size of the population).
Under mild technical conditions (\cite{EK05}, Theorem~11.2.1), the scaled process $\bX (t) / N$ converges almost surely over finite time intervals, as $N \to \infty$, to the solution $\by (t)$ of the ordinary differential equation system
\begin{eqnarray}
{d \by \over dt} &=& \sum_{\bl \in L} \bl W_{\bl} ( \by ) . \label{deterministic_ODE}
\end{eqnarray}

For our application, we suppose that the system~(\ref{deterministic_ODE}) possesses two equilibrium points: a stable endemic equilbrium point $\by^*$ with all components strictly positive, and an unstable disease-free equilibrium point at $\by = {\bf 0}$.
We next summarise some key results from the Hamiltonian approach, in a form suited to our application.
Detailed justifications and extensions of the method may be found in the review paper~\cite{AM16} and references therein.

The Hamiltonian of the system is defined to be
\begin{eqnarray}
H ( \by , \btheta ) &=& \sum_{\bl \in L} W_{\bl} ( \by ) \left( {\rm e}^{\btheta^T \bl} - 1 \right) . \label{Hamiltonian}
\end{eqnarray}
This Hamiltonian determines the following two complementary Hamilton-Jacobi partial differential equations:
\begin{eqnarray}
H \left( \by , {\partial V \over \partial \by} \right) \;\;=\;\; 0 \mbox{ and } H \left( {\partial U \over \partial \btheta} , \btheta \right) \;\;=\;\; 0 . \label{HJE}
\end{eqnarray}
Each of these Hamilton-Jacobi equations is a way of expressing the eigenvector equation~(\ref{QSD}) while retaining only leading order terms in the limit $N \to \infty$ (see Appendix for a brief outline of the derivations).

If we can solve either of the Hamilton-Jacobi equations~(\ref{HJE}), the leading-order asymptotic behaviour of the mean time to extinction~$\tau$ is given by
\begin{eqnarray}
\lim_{N\to\infty} {\ln \tau \over N} 
&=& V ( {\bf 0} ) - V ( \by^* ) 
\;\;=\;\; U ( {\bf 0} ) - U ( \btheta^* ) ,
\label{tau_U_V}
\end{eqnarray}
where $\by^*$ is the endemic equilibrium point of the deterministic system~(\ref{deterministic_ODE}), and $\btheta^*$ is the (assumed unique)  non-zero equilibrium point of the complementary system
\begin{eqnarray}
{d\btheta \over dt} &=& - \left. {\partial H \over \partial \by} \right|_{\by = {\bf 0}} . \label{theta_star_eqn}
\end{eqnarray}
Note that system~(\ref{deterministic_ODE}) may be recovered as ${d\by \over dt} = \left. {\partial H \over \partial \btheta} \right|_{\btheta = {\bf 0}}$.

The solutions $U ( \btheta )$, $V ( \by )$ to equations~(\ref{HJE}) are related via the Legendre transform; that is,
\begin{eqnarray*}
U ( \btheta ) = \sup_{\by} \left\{ \by^T \btheta - V ( \by ) \right\} , &&
V ( \by ) = \sup_{\btheta} \left\{ \btheta^T \by - U ( \btheta ) \right\} ,
\end{eqnarray*}
see~\cite{DM14}.

When (as is usually the case) it is not possible to find an analytical solution to either of the Hamilton-Jacobi equations~(\ref{HJE}), they may be solved numerically using the method of characteristics.
That is, we write down the following $2k$-dimensional system of ordinary differential equations: 
\begin{eqnarray}
\left. \begin{array}{rcl}
\displaystyle {d\by \over dt} &=&\displaystyle  {\partial H \over \partial \btheta} \;\;=\;\; \sum_{\bl \in L} \bl W_{\bl} ( \by ) {\rm e}^{\btheta^T \bl} ,
\\
\displaystyle {d\btheta \over dt} &=& \displaystyle - \; {\partial H  \over \partial \by} \;\;=\;\; - \sum_{\bl \in L} {\partial W_{\bl} \over \partial \by} \left( {\rm e}^{\btheta^T \bl} - 1 \right) ,
\end{array} \right\} \label{equations_of_motion}
\end{eqnarray}
referred to as the `equations of motion' of the system, and apply an appropriate numerical solver to~(\ref{equations_of_motion}).
We then have $\lim_{N\to\infty} ( \ln \tau ) / N = A$, where $A$ is the `action' integral,
\begin{eqnarray}
A \;\;=\;\; \int_{-\infty}^\infty \btheta^T {d\by \over dt} \, dt \;\;=\;\; - \int_{-\infty}^\infty \by^T {d\btheta \over dt} \, dt , \label{action_integral}
\end{eqnarray}
the integral in each case being evaluated along a trajectory from $( \by^* , {\bf 0} )$ to $( {\bf 0} , \btheta^* )$.
Note that $A = V ( {\bf 0} ) - V ( \by^* ) = U ( {\bf 0} ) - U ( \btheta^* )$.

Having set out the general Hamiltonian approach, we will now apply this technique to the infection model described in section~2 above.

\section{Asymptotic persistence time formulae}
\label{formulae}
Recall the infection model $\left\{ \bI (t) : t \ge 0 \right\}$ described in section~2, with transition rates given in table~1.
In the large population limit, the scaled infection process $\bI(t)/N$ converges almost surely, over finite time intervals, to the deterministic process $\by(t)$ satisfying the system of ordinary differential equations~(\ref{deterministic_ODE}); that is,
\begin{eqnarray}
{dy_i \over dt} &=& \beta \left( \sum_{j=1}^k \lambda_j y_j \right) \mu_i ( f_i - y_i ) - \gamma y_i \mbox{ for } i=1,2,\ldots,k. \label{ODE_system}
\end{eqnarray}

For $R_0 > 1$ there is a unique non-zero equilibrium point $\by^*$ of the system~(\ref{ODE_system}), and it is globally asymptotically stable~\cite{LY76}.
This endemic equilibrium point $\by^*$ is given by~\cite{N80} 
\begin{eqnarray}
y_i^* &=& {\mu_i f_i D ( \blambda,\bmu ) \over 1 + \mu_i D ( \blambda,\bmu ) } \mbox{ for } i=1,2,\ldots,k , \label{ystar}
\end{eqnarray}
where $D ( \blambda,\bmu )$ is the unique positive solution of 
\begin{eqnarray}
{\beta \over \gamma} \sum_{j=1}^k {\mu_j f_j \lambda_j \over 1 + \mu_j D ( \blambda,\bmu )} &=& 1 . \label{D_equation}
\end{eqnarray}

The Hamiltonian~(\ref{Hamiltonian}) corresponding to the process $\bI (t)$ is
\begin{eqnarray*}
H ( \by , \btheta ) &=& \beta \left( \sum_{j=1}^k \lambda_j y_j \right) \left( \sum_{i=1}^k \mu_i ( f_i - y_i ) \left( {\rm e}^{\theta_i} - 1 \right)   \right)  
+ \gamma \sum_{i=1}^k y_i \left( {\rm e}^{-\theta_i} - 1 \right) . 
\end{eqnarray*}

The corresponding equations of motion~(\ref{equations_of_motion}) are, for $i=1,2,\ldots,k$,
\begin{eqnarray}
{dy_i \over dt} &=& \beta \left( \sum_{j=1}^k \lambda_j y_j \right) \mu_i ( f_i - y_i ) {\rm e}^{\theta_i} 
- \gamma y_i {\rm e}^{-\theta_i} , \label{eqn_motion_y}
\\
{d\theta_i \over dt} &=& 
- \beta \lambda_i \sum_{j=1}^k 
\mu_j ( f_j - y_j ) 
\left( {\rm e}^{\theta_j} - 1 \right)  
+ \beta \left( \sum_{j=1}^k  
\lambda_j y_j \right) \mu_i \left( {\rm e}^{\theta_i} - 1 \right)  
- \gamma \left( {\rm e}^{-\theta_i} - 1 \right) .
\label{eqn_motion_theta}
\end{eqnarray}

The non-zero equilibrium point $\btheta^*$ given by~(\ref{theta_star_eqn}) satisfies
\begin{eqnarray}
\beta \lambda_i \sum_{j=1}^k 
\mu_j f_j \left( {\rm e}^{\theta_j^*} - 1 \right)  
+ \gamma \left( {\rm e}^{-\theta_i^*} - 1 \right) 
&=& 0 \mbox{ for } i=1,2,\ldots,k. \label{equilm1}
\end{eqnarray}

Setting $B = (\beta/\gamma) \sum_j f_j \mu_j \left( 1 -  {\rm e}^{\theta_j^*} \right)$, then~(\ref{equilm1}) implies that
\begin{eqnarray*}
{\rm e}^{-\theta_i^*} &=& 1 + \lambda_i B \mbox{ for } i=1,2,\ldots,k.
\end{eqnarray*}
Substituting back into equation~(\ref{equilm1}), we find that either $B=0$ (corresponding to $\btheta = {\bf 0}$) or $B = D ( \bmu , \blambda )$.
The elements of $\btheta^*$ are thus 
\begin{eqnarray*}
\theta_i^* &=& - \ln \left( 1 + \lambda_i D ( \bmu , \blambda )  \right) \mbox{ for } i=1,2,\ldots,k .
\end{eqnarray*}

So far, we have allowed for heterogeneities in both infectivity and susceptibility simultaneously.
If we restrict to only one type of heterogeneity, then it becomes possible to find an explicit formula for the action~$A$.
Our main result is the following.

{\bf Theorem 1.}
Consider the heterogeneous SIS infection model defined in section~2, with transition rates given in table~1, and suppose $R_0 > 1$.
Recall that $\tau$ denotes the mean time from quasi-stationarity to disease extinction, and that $D ( \blambda , \bmu )$ is defined to be the unique positive solution of equation~(\ref{D_equation}).
\begin{enumerate}
\item[(i)] If heterogeneity is in infectivity alone ($\bmu = {\bf 1}$), then
\begin{eqnarray}
\lim_{N \to \infty} {\ln \tau \over N} &=&
\sum_{i=1}^k f_i \ln \left( 1 + \lambda_i D ( {\bf 1} , \blambda ) \right) - {\gamma \over \beta} D ( {\bf 1} , \blambda ) .
\label{Action}
\end{eqnarray}
\item[(i)] If heterogeneity is in susceptibility alone ($\blambda = {\bf 1}$), then
\begin{eqnarray*}
\lim_{N \to \infty} {\ln \tau \over N} &=&
\sum_{i=1}^k f_i \ln \left( 1 + \mu_i D ( {\bf 1} , \bmu ) \right) - {\gamma \over \beta} D ( {\bf 1} , \bmu ) .
\end{eqnarray*}
\end{enumerate}
(Note: under the network interpretation, the assumption $\bmu = {\bf 1}$ corresponds to every individual having the same in-degree, whereas $\blambda = {\bf 1}$ corresponds to every individual having the same out-degree.)

{\em Proof.}
\begin{enumerate}
\item[(i)] Suppose that $\mu_i = 1$ for all $i$, and consider a trajectory $\left\{ \btheta(z) : 0 \le z \le D ( {\bf 1} , \blambda ) \right\}$ along which $\theta_i = - \ln \left( 1 + \lambda_i z \right)$ for $i=1,2,\ldots,k$.
Along such a trajectory, the Hamiltonian simplifies to 
\begin{eqnarray*}
H ( \by , \btheta ) &=& \beta \left( \sum_{j=1}^k \lambda_j y_j \right) \left( \sum_{i=1}^k ( f_i - y_i ) \left( {1 \over 1 + \lambda_i z} - 1 \right) \right)
+ \gamma \sum_{i=1}^k y_i \lambda_i z \\
&=& \gamma z \left( \sum_{j=1}^k \lambda_j y_j \right) \left( 1 - {\beta \over  \gamma} \sum_{i=1}^k {f_i \lambda_i \over 1 + \lambda_i z} 
+ {\beta \over \gamma} \sum_{i=1}^k   {y_i \lambda_i \over 1 + \lambda_i z}\right) .
\end{eqnarray*}
Since $z>0$ and $\sum_j \lambda_j y_j >0$ (except at endpoints of the trajectory) the Hamilton-Jacobi equation $H \left( {\partial U \over \partial \btheta} , \btheta \right) = 0$ reduces to 
\begin{eqnarray}
 1 - {\beta \over \gamma} \sum_{i=1}^k {f_i \lambda_i \over 1 + \lambda_i z} 
+ {\beta \over \gamma} \sum_{i=1}^k   {\lambda_i \over 1 + \lambda_i z} {\partial U \over \partial \theta_i} &=& 0 . \label{factorised}
\end{eqnarray}

Now along the trajectory under consideration, we have
\begin{eqnarray*}
{dU \over dz} &=& - \sum_{i=1}^k {\lambda_i \over 1 + \lambda_i z} {\partial U \over \partial \theta_i} 
\end{eqnarray*}
and so equation~(\ref{factorised}) becomes
\begin{eqnarray*}
{dU \over dz} &=& {\gamma \over \beta} - \sum_{i=1}^k {f_i \lambda_i \over 1 + \lambda_i z} \\
\Rightarrow \quad U( \btheta(z)) - U(\btheta(0)) &=& \left( {\gamma \over \beta} \right) z - \sum_{i=1}^k \int_0^z {f_i \lambda_i \over 1 + \lambda_i x} \, dx \\
\Rightarrow \quad U ( \btheta ) - U ( {\bf 0} ) &=& \left( {\gamma \over \beta} \right) z - \sum_{i=1}^k f_i \ln \left( 1 + \lambda_i z \right) \mbox{ along } \theta_i = - \ln ( 1 + \lambda_i z ) .
\end{eqnarray*}

The action is therefore given by
\begin{eqnarray*}
A \;\;=\;\; U ( {\bf 0} ) - U ( \btheta^* ) &=& \sum_{i=1}^k f_i \ln \left( 1 + \lambda_i D ( {\bf 1} , \blambda ) \right) -  {\gamma  \over \beta} D ( {\bf 1} , \blambda ) , 
\end{eqnarray*}
as required.

\item[(ii)] 
For the SIS model on a finite network, in which each individual $u$ makes contact with each other individual $v$ at rate $\beta_{uv}$, it is known that, provided infectious periods are exponentially distributed, the decay parameter of the process is unchanged under transposition of the matrix of infection rates $\left\{ \beta_{uv} \right\}$.
This follows from the property of `network duality', see~\cite{WS13,HL75,H76}.
In our context, this implies that the mean time to extinction from quasi-stationarity, $\tau$, is identical if we interchange the roles of $\blambda , \bmu$.
Hence part~(ii) of the theorem follows immediately from part~(i).

We can confirm this as follows.
With $\blambda = {\bf 1}$, the Hamiltonian may be written as
\begin{eqnarray*}
H ( \by , \btheta ) 
&=& \gamma \sum_{i=1}^k \left( {\rm e}^{\theta_i} - 1 \right) \left( {\beta \over \gamma} \left( \sum_{j=1}^k y_j \right) \mu_i ( f_i - y_i ) - y_i {\rm e}^{-\theta_i} \right) .
\end{eqnarray*}

With the convention that $y \ln y = 0$ when $y=0$, take
\begin{eqnarray}
V ( \by ) &=& \sum_{i=1}^k y_i \left( 1 + \ln y_i - \ln \left( {\beta \over \gamma} \mu_i \right) \right) - \left( \sum_{i=1}^k y_i \right) \ln \left( \sum_{i=1}^k y_i \right) \nonumber \\
&& {} + \sum_{i=1}^k (f_i - y_i) \ln (f_i - y_i)  . \label{V_susc}
\end{eqnarray}
Then
\begin{eqnarray*}
{\partial V \over \partial y_i} &=& 
\ln \left( {y_i \over  {\beta \over \gamma} \mu_i 
 (f_i - y_i) 
\left( \sum_j y_j \right)} \right) \mbox{ for } i=1,2,\ldots,k,
\end{eqnarray*}
and so
\begin{eqnarray*}
H \left( \by , {\partial V \over \partial \by} \right) &=& 0 .
\end{eqnarray*}
That is, $V ( \by )$ satisfies the relevant Hamilton-Jacobi equation.
The action is then given by 
\begin{eqnarray*}
A \;\;=\;\; V ( {\bf 0} ) - V ( \by^* ) &=& \sum_{i=1}^k f_i \ln \left( 1 + \mu_i D ( {\bf 1} , \bmu ) \right) - {\gamma \over \beta} D ( {\bf 1} , \bmu ) .
\end{eqnarray*}
As expected, we recover the formula for the case of heterogeneous infectivity, but with the roles of $\blambda , \bmu$ interchanged.

Having found the solution $V(\by)$ for the case $\blambda = {\bf 1}$, we can find the corresponding function $U(\btheta)$ as the Legendre transform of $V ( \by )$.
For $i=1,2,\ldots,k$, we have
\begin{eqnarray*}
{d \over dy_i} \left( \by^T \btheta - V ( \by ) \right) &=& \theta_i - \ln \left( {y_i \over {\beta \over \gamma} \mu_i \left( f_i - y_i \right) \left( \sum_j y_j \right)} \right) ,
\end{eqnarray*}
and so any stationary point satisfies
\begin{eqnarray}
y_i  &=& {{\beta \over \gamma} \mu_i f_i \left( \sum_j y_j \right) {\rm e}^{\theta_i}
\over 
1 + {\beta \over \gamma} \mu_i \left( \sum_j y_j \right) {\rm e}^{\theta_i}} . \label{ySP}
\end{eqnarray}

For $\btheta \in \R^k$, define the function $Q ( \bmu , \btheta )$ to be the solution of
\begin{eqnarray*}
{\beta \over \gamma} \sum_{j=1}^k {\mu_j f_j \over {\rm e}^{-\theta_j} + \mu_j Q ( \bmu , \btheta )} &=& 1 .
\end{eqnarray*}
Setting $R = ( \beta / \gamma ) \sum_j y_j$ and substituting from~(\ref{ySP}) into the definition of~$R$, we find that $R = Q ( \bmu , \btheta )$, and hence the function $\by^T \btheta - V ( \by )$ has a stationary point at
\begin{eqnarray}
y_i &=& {\mu_i f_i {\rm e}^{\theta_i} Q ( \bmu , \btheta )
\over 1 + \mu_i {\rm e}^{\theta_i} Q ( \bmu , \btheta )} . \label{SP}
\end{eqnarray}

Evaluating the function $\by^T \btheta - V( \by ) $ at the point~(\ref{SP}), we find 
\begin{eqnarray}
U ( \btheta ) &=& \sum_{i=1}^k f_i \ln \left( 1 + \mu_i {\rm e}^{\theta_i} Q ( \bmu , \btheta ) \right) - {\gamma \over \beta} Q ( \bmu , \btheta ) \label{Ususc}
\end{eqnarray} 
and can easily verify that the function~(\ref{Ususc}) does indeed satisfy $H \left( {\partial U \over \partial \btheta} , \btheta \right) = 0$.
Now $Q ( \bmu , {\bf 0} ) = D ( {\bf 1} , \bmu )$ and $Q ( \bmu , \btheta^* ) = 0$, so we once again find that
\begin{eqnarray*}
A &=& U ( {\bf 0} ) - U ( \btheta^* )
\;\;=\;\;
\sum_{i=1}^k f_i \ln \left( 1 + \mu_i D ( {\bf 1} , \bmu ) \right) - {\gamma \over \beta} D ( {\bf 1} , \bmu ) .
\end{eqnarray*}
\hfill$\square$
\end{enumerate}

Although we did not actually need to find the functions $U ( \btheta ) , V ( \by )$ in order to prove theorem~1(ii), we include them because knowledge of these functions can be of assistance in generalising and extending our results.
We will demonstrate this in theorem~3 below.

Figure~1 illustrates theorem~1 in the case of $k=2$ groups with heterogeneity in infectivity (the graph for the corresponding case with heterogeneity in susceptibility is identical, by network duality).
The exact value of $( \ln \tau ) / N$ is computed from equation~(\ref{QSD}) for total population sizes $N=100, 150, \ldots, 650$.
The action~$A$ is computed from equation~(\ref{Action}).
For comparison, we also show the action~$A_0 = (1/R_0) - 1 + \ln R_0$ computed for the homogeneous population SIS model with the same value for~$R_0$.
We see that formula~(\ref{Action}) gives a good approximation to $( \ln \tau ) / N$ for population sizes from around $N=300$ upwards.
We can also see that if we were to ignore heterogeneity and use the homogeneous population result, we would drastically over-estimate the persistence time of infection.
We demonstrate this point in theorem~2 below, as well as comparing heterogeneous populations of greater or lesser degrees of  heterogeneity.
\begin{figure}
\includegraphics[width=14cm]{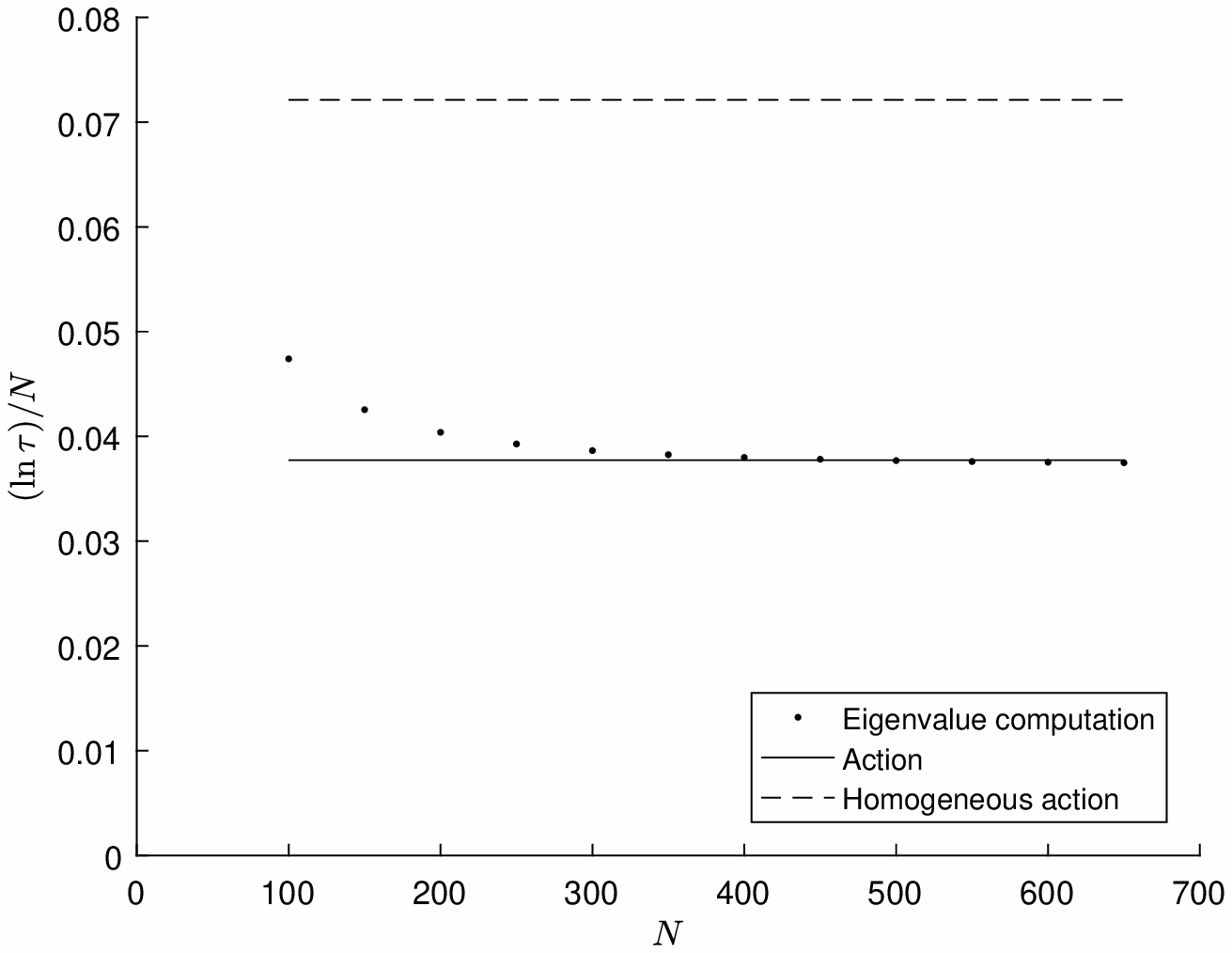}
\caption{Values of $( \ln \tau ) / N$ and asymptotic formulae plotted against population size~$N$.
Fixed parameter values $k=2$, $\bm{f} =(0.5,0.5)$, $\blambda = {2 \over 51} (50,1)$, $\bmu = (1,1)$, $R_0 = 1.5$, $\gamma = 1$.
For these parameter values, $D ( \blambda , \bmu ) = 0.5$, $\by^* = (1/6 , 1/6)$, $D ( \bmu , \blambda ) = 0.2625$, $\btheta^*=(-0.4152   -0.0102)$, , with action~$A \approx 0.0377$ and corresponding action for the homogeneous case $A_0 \approx 0.0721$.
The dots, labelled `eigenvalue computation', are the true values of $( \ln \tau ) / N$ computed from equation~(\ref{QSD}); the action~$A$ is computed from equation~(\ref{Action}); the homogeneous action is computed as $A_0 = (1/R_0) - 1 + \ln R_0$.}
\end{figure}

\section{The effect of increasing heterogeneity}
\label{effect}
Using the formulae of the previous section, we are now in a position to investigate the effect of increasing heterogeneity upon the persistence time of infection.
First, in order to compare different levels of heterogeneity, we recall the definition of majorization~\cite{MOA11}.
For any $\bx \in (\R^+)^k$, denote by $x_{[1]} \ge x_{[2]} \ge \cdots \ge x_{([k]}$ the ordered components of $\bx$.
Then for $\bx^{(1)} , \bx^{(2)} \in (\R^+)^k$, we say $\bx^{(1)}$ is majorized by $\bx^{(2)}$, denoted $\bx^{(1)} \prec \bx^{(2)}$, if $\sum_{i=1}^k x^{(1)}_i = \sum_{i=1}^k x^{(2)}_i$ and $\sum_{i=1}^j x^{(1)}_{[i]} \le \sum_{i=1}^j x^{(2)}_{[i]}$ for $j=1,2,\ldots,k-1$.
An equivalent definition is that $\sum_{i=1}^k \phi \left( x^{(1)}_i \right) \le \sum_{i=1}^k \phi \left( x^{(2)}_i \right)$ for all convex functions~$\phi(\cdot)$.
Intuitively, $\bx^{(2)}$ is `more heterogeneous' than $\bx^{(1)}$.
More generally, given a probability vector (with components summing to~1) $\bp \in ( \R^+ )^k$, then $\bx^{(1)}$ is $\bp$-majorized by $\bx^{(2)}$, written $\bx^{(1)} \prec_{\bp} \bx^{(2)}$, if there exists a permutation $\sigma$ such that 
$x^{(1)}_{\sigma(1)} \ge x^{(1)}_{\sigma(2)} \ge \cdots \ge x^{(1)}_{\sigma(k)}$ 
and $x^{(2)}_{\sigma(1)} \ge x^{(2)}_{\sigma(2)} \ge \cdots \ge x^{(2)}_{\sigma(k)}$ 
with 
$\sum_{i=1}^k p_i x^{(1)}_i = \sum_{i=1}^k p_i x^{(2)}_i$ and $\sum_{i=1}^j p_{\sigma(i)} x^{(1)}_{\sigma(i)} \le \sum_{i=1}^j p_{\sigma(i)} x^{(2)}_{\sigma(i)}$ for $j=1,2,\ldots,k-1$.

{\bf Theorem 2.}
Consider two populations, with $\beta^{(1)} = \beta^{(2)} = \beta$, $\gamma^{(1)} = \gamma^{(2)} = \gamma$, and each having the same group structure $\bm{f}^{(1)} = \bm{f}^{(2)} = \bm{f}$, where we use
superscripts $(1), (2)$ to denote the population under consideration.
Recall that~$\tau$ denotes the mean time from quasi-stationarity to disease extinction.
\begin{itemize}
\item[(i)] With heterogeneity in infectivity alone, 
\begin{eqnarray*}
\blambda^{(1)} \prec_{\bm{f}} \blambda^{(2)} &\Rightarrow& \lim_{N \to \infty} {\ln \tau^{(1)} \over N} \ge
\lim_{N \to \infty} {\ln \tau^{(2)} \over N} .
\end{eqnarray*}
\item[(ii)] With heterogeneity in susceptibility alone, 
\begin{eqnarray*}
\bmu^{(1)} \prec_{\bm{f}} \bmu^{(2)} &\Rightarrow& \lim_{N \to \infty} {\ln \tau^{(1)} \over N} \ge
\lim_{N \to \infty} {\ln \tau^{(2)} \over N} .
\end{eqnarray*}
\end{itemize}
In particular, provided heterogeneity is in either infectivity or susceptibility but not both, then $\lim_N \to \infty \left( \ln \tau \right) / N$ is maximised in the homogeneous case.

{\em Proof.}
Consider the case of heterogeneity in infectivity, and suppose that $\blambda^{(1)} \prec_{\bm{f}} \blambda^{(2)}$.
The function $h(x) = x / \left( 1 + x D ( {\bf 1} , \blambda^{(1)} ) \right)$ is concave for $x>0$, and so applying proposition~14.A.3 of~\cite{MOA11}, 
\begin{eqnarray*}
\sum_{i=1}^k {f_i \lambda_i^{(2)} \over 1 + \lambda_i^{(2)} D ( {\bf 1} , \blambda^{(1)} )}
&\le&
\sum_{i=1}^k {f_i \lambda_i^{(1)} \over 1 + \lambda_i^{(1)} D ( {\bf 1} , \blambda^{(1)} )}
\;\;=\;\; {\gamma \over \beta} ,
\end{eqnarray*}
the final equality coming from the definition~(\ref{D_equation}) of $D ( \bmu , \blambda )$ for population~1.
The expression $\sum_i f_i \lambda_i^{(2)} \left/ \left( 1 + \lambda_i^{(2)} z \right) \right.$ is a decreasing function of $z$, and so from equation~(\ref{D_equation}) for population~2 it follows that $D \left( {\bf 1} , \blambda^{(1)} \right) \ge D \left( {\bf 1} , \blambda^{(2)} \right)$.

Now define the function $\psi(z) = \sum_i f_i \ln \left( 1 + \lambda_i^{(1)} z \right) - (\gamma z / \beta)$.
Then 
\begin{eqnarray*}
{d\psi \over dz} &=& \sum_i {f_i \lambda_i^{(1)} \over 1 + \lambda_i^{(1)} z} - {\gamma \over \beta} ,
\end{eqnarray*}
so that $\psi^\prime (z) > 0$ for $0 < z < D \left( {\bf 1} , \blambda^{(1)} \right)$, and hence $\psi \left( D \left( {\bf 1} , \blambda^{(2)} \right) \right) \le \psi \left( D \left( {\bf 1} , \blambda^{(1)} \right) \right) = A^{(1)}$.
That is,
\begin{eqnarray}
\sum_{i=1}^k f_i \ln \left( 1 + \lambda_i^{(1)} D \left( {\bf 1} , \blambda^{(2)} \right) \right) - {\gamma D \left( {\bf 1} , \blambda^{(2)} \right) \over \beta}
&\le& A^{(1)} . \label{ineq1}
\end{eqnarray}

The function $g(x) = \ln ( 1 + D \left( {\bf 1} , \blambda^{(2)} \right) x )$ is concave for $x>0$, and so again applying proposition~14.A.3 of~\cite{MOA11},
\begin{eqnarray*}
\sum_{i=1}^k f_i \ln \left( 1 + \lambda_i^{(2)} D \left( {\bf 1} , \blambda^{(2)} \right) \right) 
&\le&
\sum_{i=1}^k f_i \ln \left( 1 + \lambda_i^{(1)} D \left( {\bf 1} , \blambda^{(2)} \right) \right) .
\end{eqnarray*}
Combining with~(\ref{ineq1}) yields
\begin{eqnarray*}
A^{(2)} &\le& A^{(1)} ,
\end{eqnarray*}
and the result follows.

Part~(ii) of the theorem follows immediately by interchanging the roles of $\blambda, \bmu$.
\hfill$\square$

\begin{figure}
\includegraphics[width=14cm]{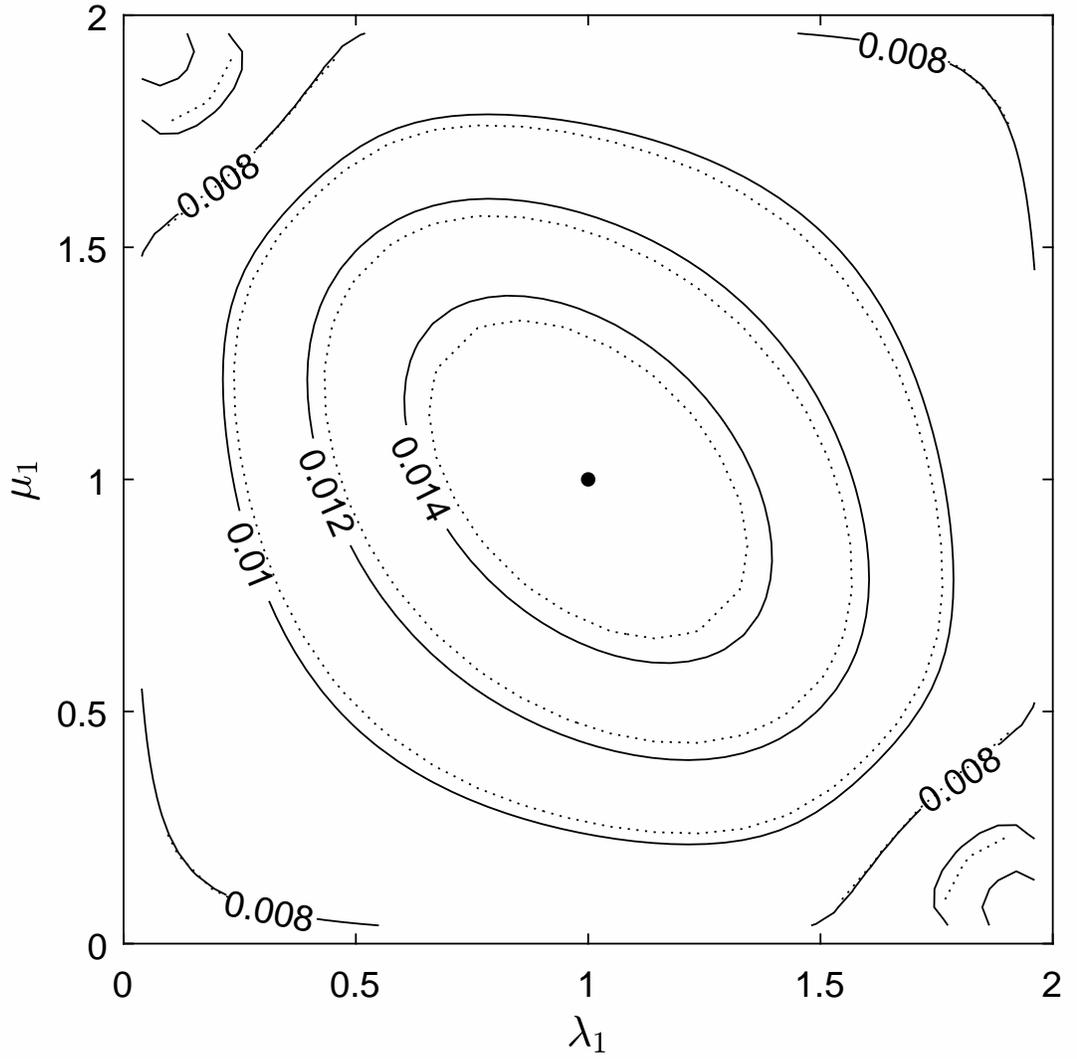}
\caption{Contour plot showing the action~$A$ (solid contours) as a function of $\lambda_1 , \mu_1$.
Fixed parameter values $k=2$, $\bm{f}=(0.5,0.5)$, $R_0 = 1.2$.
The action is maximised at $(\lambda_1 , \mu_1 ) = (1,1)$, with value $A_0 = (1/R_0) - 1 + \ln R_0 \approx 0.0157$.
For comparison, dotted contours show a finite-population approximation --- see main text for details.
}
\label{contour}
\end{figure}

Figure~\ref{contour} illustrates theorem~2, as well as showing the effect of allowing heterogeneity in both infectivity and susceptibility simultaneously, for the case of $k=2$ equal-sized groups ($f_1 = f_2 = 0.5$).
The constraints on the elements of $\blambda , \bmu$ in this case reduce to $\lambda_1 + \lambda_2 = \mu_1 + \mu_2 = 2$, and so we plot the action as a function of $( \lambda_1 , \mu_1 ) \in (0,2)^2$.
We choose to keep $R_0$ fixed, with the value of $\beta$ being varied in order to achieve this.
With both heterogeneities present, we have no explicit formula for the action~$A$, and instead compute it by first solving the equations of motion~(\ref{eqn_motion_y},\ref{eqn_motion_theta}) numerically using the Matlab {\tt bvp4c} command, and then integrating the numerical solution along the trajectory,  equation~(\ref{action_integral}).
The solid contours in figure~\ref{contour} show the action values~$A$ computed in this way.
Note that the transformation $(\lambda_1 , \mu_1 ) \to ( 2-\lambda_1 , 2-\mu_1 )$ here amounts to simply re-labelling the groups, so that figure~\ref{contour} is invariant under a rotation of half a turn around the point~$(1,1)$; also, we know from network duality that the action is unchanged under the transformation $(\lambda_1 , \mu_1) \to (\mu_1 , \lambda_1)$, so that figure~\ref{contour} is invariant under reflection in the line $\lambda_1 = \mu_1$.

For comparison, the dotted contours in figure~\ref{contour} were computed by solving the eigenvalue equation~(\ref{QSD}) numerically for $N=400$ and $N=500$, and assuming (without proof) that asymptotic formula~(\ref{general_asymptotic_formula}) is valid for our model. 
Denoting by $\tau_N$ the mean time from quasi-stationarity to disease extinction in a population of size~$N$, formula~(\ref{general_asymptotic_formula}) implies that the action~$A$ may be approximated by
\begin{eqnarray}
\left. \left( \ln \left( \tau_{500} \sqrt{500} \right) - \ln \left( \tau_{400} \sqrt{400} \right) \right) \right/ 100 ,
\label{dotted_contours}
\end{eqnarray}
and the dotted contours show computed values of formula~(\ref{dotted_contours}).
The fact that the dotted contours closely follow the solid contours provides some confirmation both that the action~$A$ gives a good approximation to $(\ln \tau ) / N$ for population sizes above $N=400$, and that  formula~(\ref{general_asymptotic_formula}) does indeed apply to our model.

We see from figure~\ref{contour} that the action decreases as we move away from the point $( \lambda_1 , \mu_1 ) = (1,1)$, not only along the lines $\lambda_1 = 1$ and $\mu_1 = 1$, as ensured by theorem~2, but in any direction.
That is, heterogeneity in infectivity, or susceptibility, or any combination of the two, reduces the value of  $\lim_{N \to \infty} ( \ln \tau ) / N$ compared to the homogeneous case.
We discuss this further in section~\ref{discussion} below.

\section{Generalising the infectious period distribution}
\label{generalising}
So far, we have made the conventional assumption that individuals' infectious periods are exponentially distributed.
This is purely a mathematical convenience, not motivated by biological realism.
Realism can be greatly improved by allowing infectious periods to follow an Erlang distribution, using the `method of stages'.
That is, when an individual becomes infected, it passes through $s$ infectious stages, remaining in each stage for an exponentially distributed time of mean $(s\gamma)^{-1}$, before returning to susceptibility.
As before, we denote by $N_j$ the (constant) number of individuals in group~$j$, and by $N = N_1 + N_2 + \ldots + N_k$ the total population size.
Denoting by $I_{jv} (t)$ the number of group~$j$ individuals in infectious stage~$v$ at time $t$, then $\left\{ I_{jv} (t) : j=1,2,\ldots,k, \ v=1,2,\ldots,s,\ t \ge 0 \right\}$ is a continuous-time Markov chain with transition rates given in table~\ref{Erlang_rates}.
The number of susceptible individuals in group~$j$ is $S_j (t) = N_j - \sum_{v=1}^s I_{jv} (t)$.

\renewcommand{\arraystretch}{2}
\begin{table}
\hspace*{-3cm}
\begin{tabular}{lll} \hline
Event & State transition & Transition rate \\ \hline
Infection in group $j$ & $I_{j1} \to I_{j1} + 1$ & ${\beta \over N} \left( \sum_{m=1}^k \lambda_m \sum_{v=1}^s I_{mv} \right) \mu_j \left( N_j - \sum_{v=1}^s I_{jv} \right)$ \\
Transition to next infectious stage & $\left( I_{j,v-1} , I_{jv} \right) \to \left( I_{j,v-1} - 1 , I_{jv} + 1 \right)$ & $s \gamma I_{j,v-1}$ for $v=2,3,\ldots,s$\\
Recovery in group $j$ & $I_{js} \to I_{js} - 1$ & $s \gamma I_{js}$ \\
\hline
\end{tabular}
\caption{Transition rates for the $k$-group, $s$-stage SIS model.}\label{Erlang_rates}
\end{table}
\renewcommand{\arraystretch}{1}

Writing $\by = \left\{ y_{iv} : i=1,2,\ldots,k,\ v=1,2,\ldots,s \right\}$ and $\btheta = \left\{ \theta_{iv} : i=1,2,\ldots,k,\ v=1,2,\ldots,s \right\}$, the corresponding Hamiltonian is 
\begin{eqnarray*}
H ( \by , \btheta ) &=& \beta
\left( \sum_{j=1}^k \lambda_j  \sum_{v=1}^s y_{jv} \right) 
\left( \sum_{i=1}^k 
 \mu_i \left( f_i - \sum_{v=1}^s y_{iv} \right) \right)  
\left( {\rm e}^{\theta_{i1}} - 1 \right) \\
&& {}
+ s \gamma \sum_{i=1}^k \sum_{v=1}^{s-1}  y_{iv} \left( {\rm e}^{-\theta_{iv}+\theta_{i,v+1}} - 1 \right)
+ s \gamma \sum_{i=1}^k y_{is} \left( {\rm e}^{-\theta_{is}} - 1 \right) .
\end{eqnarray*}

It is immediate from equation~(17) of~\cite{C15} that the deterministic endemic equilibrium point is given by
\begin{eqnarray*}
y_{iv}^* &=& {y_i^* \over s} \mbox{ for } i=1,2,\ldots,k,\ v=1,2,\ldots,s,
\end{eqnarray*}
where $y_i^*$ is the solution~(\ref{ystar}) for the model with exponentially distributed infectious periods ($s=1$).

It is straightforward to show that the elements of $\btheta^*$ are given by
\begin{eqnarray*}
\theta_{iv}^* &=& - (s+1-v) \ln \left( 1 + \lambda_i D_s ( \bmu , \blambda ) \right)
\end{eqnarray*}
where $D_s ( \bmu , \blambda )$ is the solution of
\begin{eqnarray}
{\beta \over s \gamma} \sum_{j=1}^k \mu_j f_j \lambda_j \sum_{v=1}^s \left( {1 \over 1 + \lambda_j D_s ( \bmu , \blambda )} \right)^v &=& 1 .
\end{eqnarray}

For the SIS model with Erlang-distributed infectious periods in a homogeneous population ($k=1$), the solution $U ( \btheta )$ to the relevant Hamilton-Jacobi equation was found in~\cite{CT17} to be
\begin{eqnarray*}
U ( \btheta ) &=& \ln \left( \sum_{v=1}^s {\rm e}^{\theta_v} \right) + {\gamma \over \beta} \left( {s \over \sum_{v=1}^s {\rm e}^{\theta_v}} \right) .
\end{eqnarray*}

Taking the Legendre transform, we find that $V ( \by )$ for this homogeneous-population model is given by
\begin{eqnarray}
V ( \by ) &=& \sup_{\btheta} \left\{ \btheta^T \by - U ( \btheta ) \right\} \nonumber \\
&=& \sum_{v=1}^s y_v \left( 1 + \ln y_v - \ln \left( {\beta \over s \gamma} \right) \right)
- \left( \sum_{v=1}^s y_v \right) \ln \left( \sum_{v=1}^s y_v \right) \nonumber \\
&& {}
+ \left( 1 - \sum_{v=1}^s y_v \right) \ln \left( 1 - \sum_{v=1}^s y_v \right) . \label{V_Erlang}
\end{eqnarray}

Comparing solution~(\ref{V_Erlang}) for the SIS model with Erlang-distributed infectious periods in a homogeneous population and solution~(\ref{V_susc}) for the SIS model with exponentially distributed infectious periods and heterogeneous susceptibilities, one may now guess the form of the solution $V ( \by )$ for the SIS model with Erlang-distributed infectious periods and heterogeneous susceptibilities, and verify that the relevant Hamilton-Jacobi equation is indeed satisfied.
The solution is thus found to be
\begin{eqnarray}
V ( \by ) &=& \sum_{i=1}^k \sum_{v=1}^s y_{iv} \left( 1 + \ln y_{iv} - \ln \left( {\beta \over s\gamma} \mu_i \right) \right)
- \left( \sum_{i=1}^k \sum_{v=1}^s y_{iv} \right) \ln \left( \sum_{i=1}^k \sum_{v=1}^s y_{iv} \right) \nonumber
\\ && {}
+ \sum_{i=1}^k \left( f_i - \sum_{v=1}^s y_{iv} \right) \ln \left( f_i - \sum_{v=1}^s y_{iv} \right) . \label{V_Erlang_susc}
\end{eqnarray}

Taking the Legendre transform, we find
\begin{eqnarray*}
U ( \btheta ) &=& 
\sum_{i=1}^k f_i \ln \left( 1 + \mu_i \left( \sum_{v=1}^s {\rm e}^{\theta_{iv}} \right) Q_s ( \bmu , \btheta ) \right) - {\gamma s \over \beta} Q_s ( \bmu , \btheta )
\end{eqnarray*}
where $Q_s ( \bmu, \btheta )$ is the solution of
\begin{eqnarray*}
{\beta \over \gamma s} \sum_{j=1}^k {\mu_j f_j \left( \sum_{v=1}^s {\rm e}^{\theta_{jv}} \right)
\over 
1 + \mu_j \left( \sum_{v=1}^s {\rm e}^{\theta_{jv}} \right) Q_s ( \bmu , \btheta )}
&=& 1 .
\end{eqnarray*}

The action~$A$ in this case is thus
\begin{eqnarray*}
A &=& V ( {\bf 0} ) - V ( \by^* ) 
\;\;=\;\; U ( {\bf 0} ) - U ( \btheta^* ) \;\;=\;\; \sum_i f_i \ln \left( 1 + \mu_i D ( {\bf 1} , \bmu ) \right) - {\gamma \over \beta} D ( {\bf 1} , \bmu ) ,
\end{eqnarray*}
as before, and the following result is immediate.

{\bf Theorem 3.}
Theorem~1(ii) remains valid if infectious periods are allowed to follow an Erlang, rather than exponential, distribution.
Consequently, theorem~2(ii) likewise remains valid with Erlang-distributed infectious periods.
\hfill$\square$

\begin{figure}
\includegraphics[width=14cm]{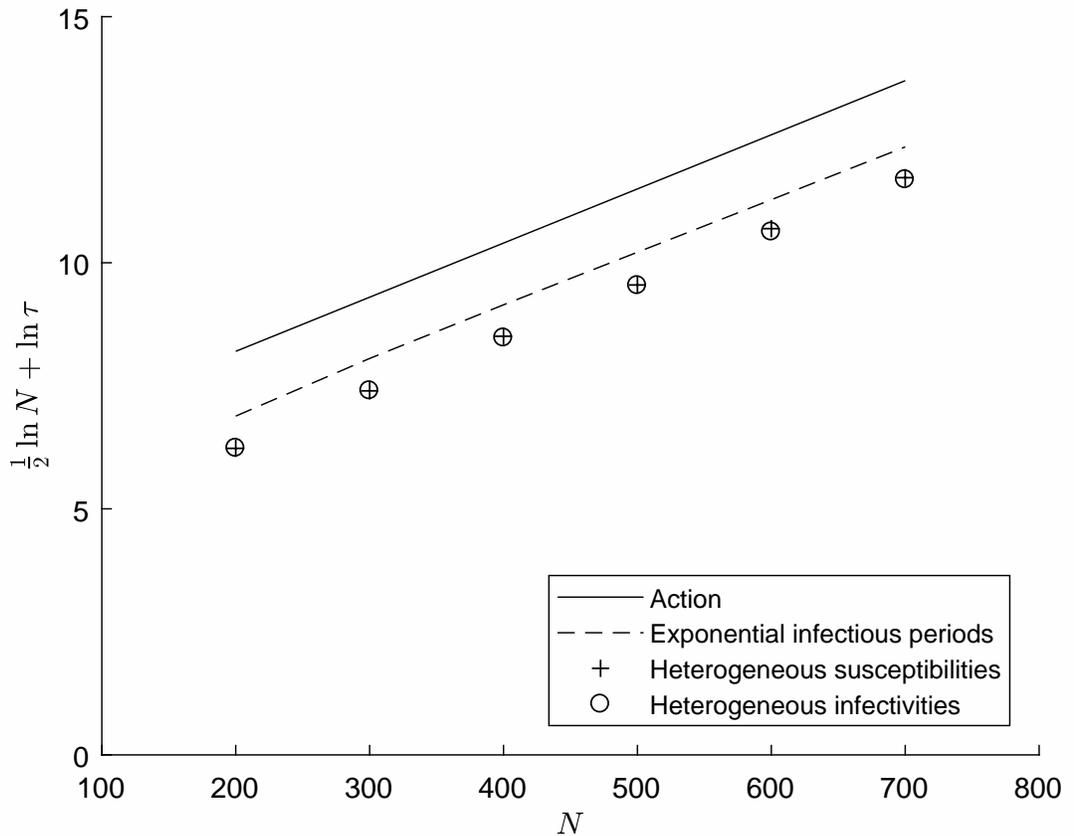}
\caption{The effect of the infectious period distribution upon the mean persistence time of infection~$\tau$.
Fixed parameter values $k=2$, $\bm{f} = (0.5,0.5)$, $R_0 = 1.2$, $\gamma = 1$.
Solid line (`Action') has gradient $A\approx 0.0110$ given by equation~(\ref{Action}) with $\blambda = (5/3,1/3)$, $\bmu = (1,1)$, intercept chosen arbitrarily; dashed line (`Exponential infectious periods')  computed from equation~(\ref{QSD}) with $\blambda = (5/3,1/3)$, $\bmu = (1,1)$; crosses (`Heterogeneous susceptibilities')  computed via simulation with $\blambda = ( 1,1 )$, $\bmu = (5/3,1/3)$ and constant infectious periods; circles (`Heterogeneous infectivities') computed via simulation with $\blambda = (5/3,1/3)$, $\bmu = (1,1)$ and constant infectious periods.}
\label{const_inf_periods}
\end{figure}

Figure~\ref{const_inf_periods} illustrates the effect of the infectious period distribution in the case of $k=2$ equal-sized groups ($f_1 = f_2 = 0.5$).
In constructing this figure, we have assumed (without proof) that asymptotic formula~(\ref{general_asymptotic_formula}) is valid for our model.
Consequently, we plot the function ${1 \over 2} \ln N + \ln \tau$, which according to formula~(\ref{general_asymptotic_formula}) should, as $N$ increases, approach a straight line of gradient~$A$ and intercept $\ln C$.
The dashed line, corresponding to  exponentially distributed infectious periods, was computed using the eigenvalue characterisation~(\ref{QSD}).
By network duality, this dashed line may be interpreted as corresponding to either heterogeneous infectivity ($\blambda = (5/3,1/3)$, $\bmu = (1,1)$) or heterogeneous suceptibility ($\blambda = (1,1)$, $\bmu = (5/3,1/3)$).
We used Monte Carlo simulation to estimate the mean persistence time~$\tau$ with constant (non-random) infectious periods
for the cases of heterogeneous infectivity and heterogeneous susceptibility separately.
This infectious period distribution corresponds to an Erlang distribution with $s$ stages in the limit as $s \to \infty$.
An issue that arises is that the time until extinction of infection, starting from quasi-stationarity, is exponentially distributed with mean increasing exponentially in population size, so that to simulate the process to extinction can be very time-consuming.
To get around this, we fixed times $t_0$ (the burn-in period) and $t_{max}$ such that (i)~by time $t_0$ the state of the process is approximately quasi-stationary (having started at time zero close to the re-scaled deterministic equilibrium point $N \by^*$, and conditioning upon survival to time $t_0$); and~(ii) by time $t_{max}$ a substantial proportion of all simulations have reached extinction.
We then estimated the mean time to extinction~$\tau$ using the maximum likelihood estimator.
That is, denoting by $T_1 , T_2 , \ldots , T_r$ the extinction times of those simulations that went extinct within the time window $\left(t_0 , t_{max} \right)$, and by $m$ the number of simulations that had not gone extinct by time $t_{max}$, our estimate is
\begin{eqnarray*}
\hat \tau &=& {m  \left( t_{max} - t_0 \right) + \sum_{i=1}^r ( T_i - t_0 ) \over r} .
\end{eqnarray*}
We have included in figure~3 a solid line with gradient equal to the action~$A$ computed from formula~(\ref{Action}).
Note that the intercept of this line was chosen arbitrarily, since we have no way to evaluate the constant~$C$ in formula~(\ref{general_asymptotic_formula}) for our model.

We see from figure~3 that the dashed line corresponding to exponentially distributed infectious periods does indeed appear to be a straight line of gradient~$A$, providing some confirmation both that the action~$A$ gives a reasonable approximation to $( \ln \tau ) / N$ for population sizes above $N=200$ and that formula~(\ref{general_asymptotic_formula}) is valid for our model.
With constant infectious periods, we see that heterogeneous infectivity and heterogeneous susceptibility result in almost identical estimates of~$\tau$, and that these estimates lie close to a straight line of gradient~$A$.
It appears that, as with exponentially distributed infectious periods, the value of~$\tau$ is unchanged if we interchange~$\blambda , \bmu$.
We therefore conjecture that, similarly to theorem~3 above, theorem~1(i) and theorem~2(i) remain valid with Erlang-distributed infectious periods.
The model with constant infectious periods has reduced mean persistence time~$\tau$ compared to the model with exponentially distributed infectious periods, but the difference is in the pre-factor constant~$C$ and not the leading-order constant~$A$, in line with the results of~\cite{BBN16} for the homogeneous population case.

\section{Discussion and possible extensions}
\label{discussion}
The main result of this paper, theorem~1, provides a simple explicit formula for $\lim_{N \to \infty} ( \ln \tau ) / N$, where $\tau$ is the expected time from endemicity to extinction for an SIS infection model with heterogeneity in either infectivity or susceptibility of individuals, in a population of size~$N$.
The only infection model for which such a formula has previously been available is the SIS model in a homogeneous population, either with exponentially distributed infectious periods~\cite{AD98} or with arbitrary infectious period distribution~\cite{BBN16}.
Theorem~1 thus represents a significant advance, but many open questions remain.

Firstly, for the SIS model in a homogeneous population, both of~\cite{AD98,BBN16} established asymptotic approximations for $\tau$ of the form~(\ref{general_asymptotic_formula}), with explicit formulae for the pre-factor constant~$C$.
Our result is less precise than this; we have not shown that an asymptotic formula of the form~(\ref{general_asymptotic_formula}) is valid for our model (although we conjecture, and have presented some numerical evidence, that this is the case), nor have we attempted to evaluate the pre-factor constant~$C$.
The asymptotic form~(\ref{general_asymptotic_formula}) has been shown by~\cite{AM10} to be valid (and formulae given for the constant~$C$) for general 1-dimensional processes of bounded jump size.
The technique of~\cite{AM10} is an extension of the approach employed here, retaining terms beyond the leading order in~$N$.
The analysis is considerably more intricate than the leading-order treatment we have restricted ourselves to, and it is not clear how the approach of~\cite{AM10} may be extended to multi-dimensional processes.

Secondly, our model as described in section~2 incorporates heterogeneity in both infectivity and susceptibility simultaneously, but we have only been able to provide an explicit asymptotic formula for the cases in which only one of these two heterogeneities is present.
In particular, this severely restricts the class of networks to which our results may be applied under the annealed network approximation --- we require either that every individual has the same in-degree, or that every individual has the same out-degree.
Nevertheless, for the class of directed networks to which they apply our results represent an interesting step forward, and since network models are of great ongoing interest in infection modelling, we now present our results in a form suited to the network interpretation.

As a preliminary to the statement of our results, we require the concept of convex ordering of random variables, defined as follows (\cite{SS07}, section~3.A.1).
Given two random variables $X^{(1)},X^{(2)}$, then $X^{(2)}$ is greater than $X^{(1)}$ in the sense of convex ordering, denoted $X^{(1)} \le_{cv} X^{(2)}$, if $E \left[ \phi \left( X^{(1)} \right) \right] \le E \left[ \phi \left( X^{(2)} \right) \right]$ for all convex functions $\phi(\cdot)$.
If $X^{(1)} , X^{(2)}$ take values in $\left\{ 1,2,\ldots , d_{\max} \right\}$, then an equivalent definition is that
\begin{eqnarray*}
\sum_{i=1}^j P \left( X^{(1)} \le i \right)
&\le&
\sum_{i=1}^j P \left( X^{(2)} \le i \right)
\mbox{ for } j=1,2,\ldots,d_{\max} .
\end{eqnarray*}
Note that $X^{(1)} \le_{cv} X^{(2)}$ implies $E \left[ X^{(1)} \right] = E \left[ X^{(2)} \right]$; intuitively, $X^{(2)}$ is `more variable' than $X^{(1)}$.

In comparing two populations, we  need to define our `groups' slightly differently than in section~2.
Specifically, partition the population into groups in such a way that two individuals belong to the same group if they share the same values of both $\left( d_{\mbox{in}}^{(1)} , d_{\mbox{out}}^{(1)} \right)$ and $\left( d_{\mbox{in}}^{(2)} , d_{\mbox{out}}^{(2)} \right)$; the condition $\bm{f}^{(1)} = \bm{f}^{(2)} = \bm{f}$ required by theorem~2 is thus satisfied.
(As before, superscripts $(1),(2)$ denote the population under consideration.)

{\bf Theorem 4.}
Consider an SIS infection in a population of $N$~individuals connected by an uncorrelated directed network.
Each individual has in-degree and out-degree distributed as $\left( d_{\mbox{in}}, d_{\mbox{out}} \right)$,  the degrees of distinct individuals being mutually independent, with $E \left[ d_{\mbox{in}} \right] = E \left[ d_{\mbox{out}} \right] = \mu$ and $d_{\mbox{in}} , d_{\mbox{out}} \le d_{\max}$ for some $d_{\max} \in \N$.
Infection transmits along each link from an infectious to a susceptible individual at rate $\kappa$, and when an individual becomes infected it remains so for a time of mean $1/\gamma$ before returning to the susceptible state.
Recall that $\tau$ denotes the expected time from quasi-stationarity to extinction of infection.

\begin{enumerate}
\item[(a)] Suppose that $P \left( d_{\mbox{in}} = \mu \right) = 1$, so every individual has the same in-degree, and that infectious periods are exponentially distributed.
Then
\begin{enumerate}
\item[(i)] $\displaystyle \lim_{N \to \infty} {\ln \tau \over N} \approx A_{\mbox{out}}$, with
\begin{eqnarray*}
A_{\mbox{out}} &=& \sum_{i=1}^{d_{\max}} 
P \left( d_{\mbox{out}} = i \right) 
\ln \left( 1 + i D_{\mbox{out}} \right) - {\gamma \over \kappa} D_{\mbox{out}}
\end{eqnarray*}
where $D_{\mbox{out}}$ is the solution of
\begin{eqnarray*}
{\kappa \over \gamma} \sum_{j=1}^{d_{\max}} {j P \left( d_{\mbox{out}} = j \right) \over 1 + j D_{\mbox{out}}} &=& 1 ;
\end{eqnarray*}
\item[(ii)] for two populations with $\kappa^{(1)} = \kappa^{(2)}$, $\gamma^{(1)} = \gamma^{(2)}$,
\begin{eqnarray*}
d^{(1)}_{\mbox{out}} \le_{cv} d^{(2)}_{\mbox{out}} 
&\Rightarrow&
A_{\mbox{out}}^{(1)} \ge A_{\mbox{out}}^{(2)} ;
\end{eqnarray*} 
in particular, $A_{\mbox{out}}$ is maximised when every individual has the same out-degree, $P \left( d_{\mbox{out}} = \mu \right) = 1$.
\end{enumerate}
\item[(b)] Suppose that $P \left( d_{\mbox{out}}  = \mu \right) = 1$, so every individual has the same out-degree, and that infectious periods follow an Erlang distribution.
Then
\begin{enumerate}
\item[(i)] $\displaystyle \lim_{N \to \infty} {\ln \tau \over N} \approx A_{\mbox{in}}$, with
\begin{eqnarray*}
A_{\mbox{in}} &=&
\sum_{i=1}^{d_{\max}} 
P \left( d_{\mbox{in}} = i \right) 
\ln \left( 1 + i D_{\mbox{in}} \right) - {\gamma \over \kappa} D_{\mbox{in}}
\end{eqnarray*}
where $D_{\mbox{in}}$ is the solution of 
\begin{eqnarray*}
{\kappa  \over \gamma} \sum_{j=1}^{d_{\max}} 
{j P \left( d_{\mbox{in}} = j \right) \over 1 + j D_{\mbox{in}}} &=& 1 ;
\end{eqnarray*}
\item[(ii)] for two populations with $\kappa^{(1)} = \kappa^{(2)}$, $\gamma^{(1)} = \gamma^{(2)}$, 
\begin{eqnarray*}
d^{(1)}_{\mbox{in}} \le_{cv} d_{\mbox{in}}^{(2)} 
&\Rightarrow&
A^{(1)}_{\mbox{in}} \ge A^{(2)}_{\mbox{in}} ;
\end{eqnarray*}
in particular, $A_{\mbox{in}}$ is maximised when every individual has the same in-degree, $P \left( d_{\mbox{in}} = \mu \right) = 1$.
\end{enumerate}
\end{enumerate}

Note that our asymptotic results (theorem~1) are exact for the model with transition rates given by table~\ref{k_group_SIS}, but approximate under the annealed network interpretation.

A third open issue is to allow for more general infectious period distributions in the case of heterogeneous infectivity.
We conjecture that theorem~1(i) remains valid, and hence also theorem~2(i), if infectious periods are allowed to follow an Erlang, rather than exponential, distribution.
Indeed, figure~\ref{const_inf_periods} suggests that the mean persistence time~$\tau$ is unchanged when $\blambda , \bmu$ are interchanged.
A difficulty here is that, in contrast to the case of heterogeneous susceptibilities, we have not been able to find complete solutions $U ( \btheta ) , V ( \by )$ to the relevant Hamilton-Jacobi equations even for the case of exponentially distributed infectious periods, but only, in proving theorem~1(i), to evaluate $U ( \btheta )$ along one particular trajectory.

One advantage of simple asymptotic formulae such as provided by theorem~1, as opposed to the exact formula~(\ref{QSD}), is that they provide a route to qualitative results such as theorem~2, that increasing heterogeneity reduces (at least to leading order) the expected persistence time of infection, and in particular that persistence time is maximised in a homogeneous population.
Theorem~2 establishes this ordering when heterogeneity is in either infectivity or susceptibility; figure~\ref{contour} suggests that the result remains true even when both types of heterogeneity are present.
It is interesting to compare with the results contained in section~5 of~\cite{CP13} regarding the effect of such heterogeneities upon the (large population) mean endemic prevalence level~$y^* = \sum_{i=1}^k y_i^*$.
Theorem~7(i) and theorem~10 of~\cite{CP13} show, respectively, that heterogeneous infectivity alone has no effect upon the endemic prevalence level~$y^*$, whereas heterogeneous suceptibility alone can only decrease~$y^*$, with $\bmu^{(1)} \prec_{\bm{f}} \bmu^{(2)} \Rightarrow  y^{*(1)} \ge y^{*(2)}$, corresponding to our theorem~2(ii) for persistence times.
When both types of heterogeneity are combined, theorem~8 of~\cite{CP13} shows that if infectivity and susceptibility are non-negatively correlated ($\sum_{i=1}^k \lambda_i \mu_i f_i \ge 1$) then the endemic prevalence level cannot be greater than for a homogeneous population with the same $R_0$ value. 
However, theorem~9 of~\cite{CP13} together with numerical work shown in figure~3 of~\cite{CP13} demonstrates that when infectivity and susceptibility are negatively correlated it is possible for the endemic prevalence level to be greater than in the homogeneous case (with $R_0$ values matched).
This presents an interesting contrast to our numerical results in figure~\ref{contour}, where heterogeneities were found to always decrease (to leading order) the expected persistence time, regardless of whether $\blambda , \bmu$ are positively or negatively correlated.

There is a slightly counter-intuitive aspect to the above results, in that an increase in endemic prevalence level may correspond to a decrease in expected persistence time.
This is easily resolved by observing that an increase in prevalence level may be accompanied by a corresponding increase in variability, leading to faster extinction of infection.
The effect of heterogeneities upon the variability of the quasi-stationary distribution is studied in section~6 of~\cite{CP13}, via an Ornstein-Uhlenbeck diffusion approximation that leads to a multivariate normal approximation to the quasi-stationary distribution~$\bq$.
The variability of this approximating normal distribution is then used as a proxy measure of persistence time.
This approach seems reasonable in terms of qualitative comparisons between infection models, and is common in the literature.
However, the approach is known to give a very bad numerical approximation to mean persistence time, with incorrect leading-order asymptotic behaviour, due to the failure in the lower tail of the normal approximation to the quasi-stationary distribution~\cite{DSS05,CT17}.
The methods of the current paper, in contrast, deal directly with the expected persistence time and yield correct leading-order asymptotic formulae.

\appendix
\section{Appendix}
The Hamilton-Jacobi equations~(\ref{HJE}) amount to two alternative ways of expressing the eigenvector equation~(\ref{QSD}), retaining only terms of leading order in $N$ in the limit as $N \to \infty$.
We briefly outline the derivation of each equation, referring the reader to~\cite{AM16} and references therein for full details.

For a process with transition rates of the form~(\ref{rates}), equation~(\ref{QSD}) may be written as
\begin{eqnarray}
\sum_{\bl \in L} \left( q_{\bx-\bl}  W_{\bl} \left( {\bx - \bl \over N} \right) -  q_{\bx} W_{\bl} \left( {\bx \over N} \right) \right) &=& - (\tau N)^{-1} q_{\bx} \mbox{ for } \bx \in C . \label{A1}
\end{eqnarray} 

Suppose ({\em ansatz}) that $q_{\bx} = \exp \left( - N V ( \bx / N ) + o(N) \right)$ for some function $V(\cdot)$.
Writing $\by = \bx / N$, then $V \left( \by - {\bl \over N} \right) = V ( \by ) - ( \bl^T / N )  {\partial V \over \partial \by} + o (1/N)$ and hence
\begin{eqnarray*}
q_{\bx - \bl} &=& \left( \exp \left( \bl^T {\partial V \over \partial y} \right) + o(1) \right) q_{\bx}  .
\end{eqnarray*}
Similarly, $W_{\bl} \left( \by - {\bl \over N} \right) = W_{\bl} ( \by ) + o(1/N)$.
Retaining only terms of leading order in~$N$, and assuming that $\tau$ is sufficiently large for the right hand side of equation~(\ref{A1}) to be neglected, then equation~(\ref{A1}) reduces to
\begin{eqnarray*}
\sum_{\bl \in L}  W_{\bl} ( \by )\left( \exp \left( \bl^T {\partial V \over \partial \by} \right) - 1 \right)
&=& 0 .
\end{eqnarray*}
That is, $H \left( \by , {\partial V \over \partial \by} \right) = 0$.

Next, consider the moment generating function $M ( \btheta ) = \sum_{\bx \in C}  {\rm e}^{\btheta^T \bx} q_{\bx}$.
Multiplying equation~(\ref{A1}) by ${\rm e}^{\btheta^T \bx}$ and summing over $\bx$ we find
\begin{eqnarray}
\sum_{\bx \in C} \sum_{\bl \in L} {\rm e}^{\btheta^T \bx} q_{\bx} W_{\bl} \left( { \bx \over N} \right) \left( {\rm e}^{\btheta^T \bl} - 1 \right) &=& - (\tau N)^{-1} M ( \btheta ) . \nonumber \\
\mbox{That is,} \quad
\sum_{\bl \in L} \left( {\rm e}^{\btheta^T \bl} - 1 \right) W_{\bl} \left( {1 \over N} {\partial \over \partial \btheta} \right) M ( \btheta )  &=& - (\tau N)^{-1} M ( \btheta ) . \label{A2}
\end{eqnarray}

Suppose ({\em ansatz}) that $M ( \btheta ) = \exp \left( N U ( \btheta ) + o(N) \right)$ for some function $U(\cdot)$, so that $U(\cdot)$ gives the leading-order term in the cumulant generating function of the quasi-stationary distribution~$\bq$.
Retaining only terms of leading order in~$N$ and assuming that $\tau$ is sufficiently large for the right hand side of equation~(\ref{A2}) to be neglected, then equation~(\ref{A2}) reduces to
\begin{eqnarray*}
\sum_{\bl \in L} \left( {\rm e}^{\btheta^T \bl} - 1 \right) W_{\bl} \left( {\partial U \over \partial \btheta} \right) &=& 0 .
\end{eqnarray*}
That is, $H \left( {\partial U \over \partial \btheta} , \btheta \right) = 0$.

Note that the Hamilton-Jacobi equations~(\ref{HJE}) thus have boundary conditions $U ( {\bf 0} ) = V ( \by^* ) = 0$.
However, we can add an arbitrary constant to any solution of either equation~(\ref{HJE}) and still have a solution, and any such additive constant will cancel out when we come to estimate $(\ln \tau ) / N$ using equation~(\ref{tau_U_V}).
We therefore ignore the boundary conditions and instead choose additive constants so as to keep our presented solutions for $U( \btheta ) , V ( \by )$ as simple as possible.

\end{document}